\documentclass[11pt]{article}
\usepackage{amsfonts,amssymb,amsmath,amsthm}
\usepackage[all]{xy}
\usepackage[hidelinks,pdftex]{hyperref} 
\usepackage[bottom]{footmisc} 

\title{On an example of Aspinwall, Morrison, \\ and Szendr\H oi}
\author{Nicolas Addington and Benjamin Tighe}
\date{}

\newcommand \C {\mathbb C}
\newcommand \F {\mathbb F}
\newcommand \R {\mathbb R}
\renewcommand \O {\mathcal O}
\renewcommand \P {\mathbb P}
\newcommand \Q {\mathbb Q}
\newcommand \Z {\mathbb Z}

\renewcommand \phi {\varphi}

\DeclareMathOperator \Aut {Aut}
\DeclareMathOperator \Hom {Hom}
\DeclareMathOperator \rank {rank}
\DeclareMathOperator \coker {coker}
\DeclareMathOperator \End {End}
\DeclareMathOperator \Hdg {Hdg}
\DeclareMathOperator \Hilb {Hilb}
\DeclareMathOperator \Sp {Sp}
\DeclareMathOperator \sign {sign}
\let \Im \undefined
\DeclareMathOperator \Im {Im}

\newtheorem* {main_thm*} {Main Theorem}
\newtheorem {thm} {Theorem} [section]
\newtheorem {prop} [thm] {Proposition}
\newtheorem {lem} [thm] {Lemma}
\theoremstyle {remark}
\newtheorem {rmk} [thm] {Remark}
\newtheorem* {rmks*} {Remarks}

\renewcommand \labelenumi {(\alph{enumi})}

\makeatletter
\DeclareRobustCommand{\xRightarrow}[2][]{\ext@arrow 0359\Rightarrowfill@{#1}{#2}}
\makeatother

\begin{document}

\maketitle


\vspace{-1.5em} 

\begin{abstract}
We study the cohomology of a 1-parameter family $Y_t$ of Calabi--Yau 3-folds introduced by Aspinwall and Morrison, related to the mirror quintic family.  Szendr\H oi proved that $Y_t, Y_{\xi t}, \dotsc, Y_{\xi^4 t}$, where $\xi$ is a fifth root of unity, have the same rational Hodge structure but are not isomorphic, and conjectured that they are not birational or even derived equivalent.  We confirm this by proving that their \emph{integral} Hodge structures are different, and discuss how this fits with known Torelli-type theorems and counterexamples.
\end{abstract}

\section*{Introduction}

For a pair of Calabi--Yau 3-folds\footnote{Defined as smooth complex projective 3-folds with trivial canonical bundle and $b_1 = 0$, but we allow $\pi_1 \ne 0$.} $X$ and $Y$, each of the following conditions implies the next:
\begin{enumerate}
\renewcommand \theenumi {(\roman{enumi})}
\renewcommand \labelenumi \theenumi
\item \label{iso} $X$ and $Y$ are isomorphic.
\item \label{bir} $X$ and $Y$ are birational.
\item \label{der} $D^b_\text{coh}(X) \cong D^b_\text{coh}(Y)$.
\item \label{H3Z} $H^3(X,\Z)_\text{tf} \cong H^3(Y,\Z)_\text{tf}$ as polarized Hodge structures.\footnote{We will use the subscript tf for the torsion-free part, that is, $H^3(-,\Z)/\text{torsion}$.}
\item \label{H3Q} $H^3(X,\Q) \cong H^3(Y,\Q)$ as polarized Hodge structures.
\end{enumerate}
The implication \ref{bir} $\Rightarrow$ \ref{der} is due to Bridgeland \cite{bridgeland_flops}.  The implication \ref{der} $\Rightarrow$ \ref{H3Z} is due to the first author \cite[footnote 8]{adm}, with a fuller account given by Ottem and Rennemo in \cite[Prop.~2.1]{or}.  Various compositions were known earlier: for example, \ref{der} $\Rightarrow$ \ref{H3Q} just uses textbook techniques \cite[\S5.2]{huybrechts_fm},
while \ref{bir} $\Rightarrow$ \ref{H3Q} goes back to Koll\'ar \cite[Cor.~4.12 and Rmk.~4.13(i)]{kollar_flops}.

The reverse implications, if they held, would be Torelli-type theorems, but they tend to fail.  Atiyah flops often give counterexamples to \ref{bir} $\Rightarrow$ \ref{iso}; the first counterexamples to stay within a single deformation family were given by Szendr\H oi \cite{balazs_birational}.  Counterexamples to \ref{der} $\Rightarrow$ \ref{bir} were first given by Borisov and C\u{a}ld\u{a}raru \cite{bc}; further counterexamples given by Schnell \cite{schnell}, Hosono and Takagi \cite{ht_der_eq}, and the first author and Bragg \cite{ab} also demonstrate the need to kill torsion in \ref{H3Z}, as the first author observed in \cite{br_not_inv}; a counterexample to \ref{der} $\Rightarrow$ \ref{bir} that stays within a single deformation family was given by Ottem and Rennemo \cite{or} and independently by Borisov, C\u{a}ld\u{a}raru, and Perry \cite{bcp}.

No counterexample to \ref{H3Z} $\Rightarrow$ \ref{der} is known, but a strong candidate was constructed by Aspinwall and Morrison in \cite{am}.  They studied a family of Calabi--Yau 3-folds $Y_t$ for $t \in \C$ with $t^5 \ne 1$, closely related to the mirror family of the Fermat quintic: they took the Dwork pencil
\[ Q_t = \{ z_0^5 + z_1^5 + z_2^5 + z_3^5 + z_4^5 = 5t \cdot z_0 z_1 z_2 z_3 z_4 \} \subset \C\P^4, \]
divided by an action of $G := (\Z/5 \times \Z/5) \rtimes \Z/5$, and took a crepant resolution of singularities, as we review in \S\ref{asmosz_setup} and \S\ref{1/5(1,1,3)}.  They gave a physics computation of the genus-0 and genus-1 Gromov--Witten invariants of the mirror family, building on \cite{cogp}, which suggested that $Y_t$ and $Y_{\xi t}$, where $\xi = e^{2 \pi i/5}$, should have the same Hodge structure on $H^3$, but different derived categories of coherent sheaves.

Szendr\H oi \cite{balazs} showed that the universal covers $Z_t$ and $Z_{\xi t}$ of $Y_t$ and $Y_{\xi t}$ are isomorphic, and that $H^3(Y_t,\Q) \cong H^3(Z_t,\Q)$ up to rescaling the intersection form,
so
\begin{equation} \label{Q-Hodge-isom}
H^3(Y_t, \Q) \cong H^3(Y_{\xi t}, \Q)
\end{equation}
as polarized Hodge structures.  On the other hand, he showed that if $t$ is general then $Y_t, Y_{\xi t}, \dotsc, Y_{\xi^4 t}$ are not isomorphic, and conjectured that they are not birational \cite[Conj.~0.2]{balazs} or even derived equivalent \cite[slide 21]{balazs_simons_talk}.

We initially hoped to upgrade the rational Hodge isometry \eqref{Q-Hodge-isom} to integral coefficients, then prove that $Y_t, Y_{\xi t}, \dotsc, Y_{\xi^4 t}$ are not derived equivalent using either the categorical enumerative invariants of C\u{a}ld\u{a}raru and Tu \cite{cei}, or the BCOV invariant of Fang, Lu, and Yoshikawa \cite{fly}.  Instead, we will prove:
\begin{main_thm*}
Let $Y_t$ be the Calabi--Yau 3-folds studied by Aspinwall and Morrison \cite{am} and Szendr\H oi \cite{balazs}.  If $t$ is very general, then no two of
\begin{align*}
H^3(Y_t, \Z)_\text{tf} &&
H^3(Y_{\xi t}, \Z)_\text{tf} &&
H^3(Y_{\xi^2 t}, \Z)_\text{tf} &&
H^3(Y_{\xi^3 t}, \Z)_\text{tf} &&
H^3(Y_{\xi^4 t}, \Z)_\text{tf}
\end{align*}
are isomorphic to one another, even as unpolarized Hodge structures.  In particular, $Y_t, Y_{\xi t}, \dotsc, Y_{\xi^4 t}$ are neither derived equivalent nor birational.
\end{main_thm*}

So Aspinwall and Morrison's example, rather than showing that \ref{H3Z} does not imply \ref{der}, instead shows that \ref{H3Q} does not imply \ref{H3Z} even generically.  We doubt that \ref{H3Z} implies \ref{der}, but we do not know where to look for a counterexample.

One might well have expected \ref{H3Q} to imply \ref{H3Z} generically: it is true for quintic 3-folds, and for many hypersurfaces both Calabi--Yau and otherwise, as well as curves of genus $\ge 4$, as Voisin explained in \cite[Rmk.~0.3]{voisin_schiffer}. 
In fact it is reasonable to expect whenever the dimension of the moduli space is less than that of the period domain, as Will Sawin explained to us \cite{sawin_mo}: that excludes K3 surfaces, hyperk\"ahler varieties, cubic 4-folds, Abelian varieties, and curves of genus $\le 3$, but includes most other cases.  In our case, the moduli space is 1-dimensional while the polarized period domain is 4-dimensional.

Filip \cite[Thm.~7]{filip} proved a \emph{monodromy} Torelli theorem for several families of Calabi--Yau 3-folds including the mirror family of the quintic 3-fold, which also applies to the families $Y_t$ and $Z_t$ because they give the same rational variation of Hodge structure. 
Precisely, let $B = \C \setminus \{ 1, \xi, \xi^2, \xi^3, \xi^4 \}$ be the base of these three families, let $\Gamma$ be the monodromy group -- that is, the image of $\pi_1(B)$ in $\Sp_4(\Q)$ -- and let $D$ be the local period domain; then the map $B \to \Gamma\backslash D$ is injective.\footnote{In fact Filip used a period domain that only remembers $F^2 = H^{3,0} \oplus H^{2,1} \subset H^3$, not the usual flag $F^3 \subset F^2 \subset H^3$, but injectivity of that map implies injectivity of this one.}
The monodromy group $\Gamma$ is known to be Zariski dense in $\Sp_4(\Q)$, but not finite index in $\Sp_4(\Z)$ by \cite{bt}; in fact each of the three families yields a different copy of $\Sp_4(\Z)$ in $\Sp_4(\Q)$, but all three are commensurable.  Filip's result strengthened an earlier \emph{generic} monodromy Torelli theorem due to Usui \cite{usui}.  We do not know whether a stronger Torelli theorem holds for any of the three, but we reiterate our surprise that a Torelli theorem of any description for $Y_t$ with $\Z$ coefficients would not imply at least a generic one with $\Q$ coefficients.

\begin{rmks*} \ 
\begin{enumerate}
\item If $t$ is very general then $Y_t, Y_{\xi t}, \dotsc, Y_{\xi^4 t}$ are not L-equivalent: that is, no two of them become equal in the Grothendieck ring of varieties after multiplying by a power of the affine line.  This follows from an argument of Efimov \cite[proof of Thm.~3.1]{efimov} and Proposition \ref{endomorphisms}(b), which implies that the Hodge endomorphism ring $\End_{\Hdg}(H^3(Y_t, \Z)_\text{tf}) = \Z$.

\item We wonder whether the rational Chow motives of $Y_t, Y_{\xi t}, \dotsc, Y_{\xi^4 t}$ are isomorphic.
\end{enumerate}
\end{rmks*}

The idea of the proof is that if $t$ is very general, then up to a scalar, the only Hodge isomorphism between $H^3(Y_{\xi^a t}, \Q)$ and $H^3(Y_{\xi^b t}, \Q)$ is the one already mentioned \eqref{Q-Hodge-isom}, but this isomorphism does not preserve the integral structure, as we argue in \S\ref{Hodge_endos} and \S\ref{finis_laus_deo}.  We study that integral structure in \S\ref{H3_section}, showing that $H^3(Y_t, \Z)_\text{tf}$ is naturally identified with a finite-index subgroup of the $G$-invariant part of $H^3(Q_t,\Z)$, given by the kernel of certain boundary maps in a Cartan--Leray spectral sequence involving the group cohomology of $G$ with coefficients in $H^*(Q_t,\Z)$.

To pin everything down, we need an explicit description of $H^3(Q_t,\Z)$ and the action of $\Aut(Q_t)$ on it.  In fact it is enough to understand $t=0$, the Fermat quintic.  Walcher \cite[\S5]{walcher}, building on calculations of Brunner et al.\ \cite[\S2.3.2]{bdlr}, observed that if $\alpha \in H^3(Q_0,\Z)$ is the Poincar\'e dual of the real 3-manifold $Q_0(\R)$, then the $\Aut(Q_0)$-orbit of $\alpha$ spans $H^3(Q_0,\Z)$.  We give a self-contained account of the proof in \S\ref{fermat}, and Magma code to carry out the needed calculations in the ancillary file \verb|quintic.magma|.  The idea is to compute the intersection numbers $\langle \phi^*\alpha,\, \psi^* \alpha \rangle$ for all $\phi, \psi \in \Aut(Q_0)$, check that the rank of the resulting pairing agrees with the rank of $H^3(Q_0,\Z)$, and that after dividing by the kernel it becomes a perfect pairing.  When $\phi(Q_0(\R))$ and $\psi(Q_0(\R))$ intersect transversely, computing the intersection number is a simple matter of checking orientations; more care is needed when the intersection has $S^1$ or $\R\P^2$ as a connected component, but the fact that $Q_0(\R)$ and its translates are Lagrangian submanifolds of $Q_0$ simplifies the analysis of excess normal bundles in those cases.  Experimentally, the analogous statement seems to hold for Fermat hypersurfaces of any odd degree and dimension; it would be interesting to find a conceptual reason for this.

We could instead have used Looijenga's description of the integral cohomology of Fermat varieties \cite[Cor.~2.2]{looijenga_fermat}, building on work of Pham \cite{pham}, with a correction discussed in \cite[Rmk.~2.2]{gs} and \cite[Cor.~4.2]{looijenga2}.  Or we could have used a set of Lagrangian 3-spheres constructed by Smith, Thomas, and Yau in \cite[Lem.~3.11]{sty}, building on the celebrated paper of Candelas et al.\ \cite[App.~A]{cogp}.  But both constructions are much more intricate than simply using $Q_0(\R)$.

\paragraph*{Acknowledgements.}
We thank Paul Aspinwall, Ross Casebolt, Will \linebreak 
Donovan, Dan Dugger, Robert Lipshitz, Eduard Looijenga, Ben Moonen, Evgeny Shinder, Ivan Smith, Bal\'azs Szendr\H oi, and Richard Thomas for helpful discussions and correspondence, Johannes Walcher for drawing our attention to \cite{bdlr}, \cite{mw}, and \cite{walcher}, and a careful referee for noticing that we had looked up the cohomology of the wrong non-Abelian group of order 125.  We used Claude Opus 5 for proofreading; it flagged many typos and a minor error in an earlier version of Proposition \ref{orientation_prop}.  N.A.\ was partially supported by NSF grant no.\  DMS-1902213.  B.T.\ was supported by NSF grant no.\ DMS-2039316.


\section{Integral cohomology of the Fermat quintic 3-fold} \label{fermat}

\subsection{Preliminaries}

Consider the Fermat quintic 3-fold
\[ Q_0 = \{ z_0^5 + z_1^5 + z_2^5 + z_3^5 + z_4^5 = 0 \} \subset \C\P^4. \]
The real points $Q_0(\R)$ form a smooth real 3-manifold, and a calculation with Stiefel--Whitney classes shows that it is orientable.  In fact it is homeomorphic to $\R\P^3$: the homeomorphism from $\R\P^4$ to itself given by
\[ (x_0 : x_1 : \dotsb : x_4) \mapsto (x_0^{1/5} : x_1^{1/5} : \dotsb : x_4^{1/5}) \]
takes the hyperplane $\{ x_0 + \dotsb + x_4 = 0 \}$ to $Q_0(\R)$.\footnote{This map is not smooth, but any homeomorphism between smooth 3-manifolds can be perturbed to a diffeomorphism \cite[Thm.~6.3]{munkres}.}  Choose an orientation of $Q_0(\R)$ and let $\alpha \in H^3(Q_0,\Z)$ be the Poincar\'e dual.  In this half of the paper we give a self-contained account of the following result:

\begin{thm}[Walcher {\cite[\S5]{walcher}}] \label{real_parts}
The class $\alpha$ and its translates by $\Aut(Q_0)$ generate $H^3(Q_0,\Z)$.
\end{thm}

The automorphism group of $Q_0$ contains a copy of the symmetric group $S_5$ that permutes the coordinates, and copies of $\Z/5$ that act on each coordinate by $\xi = e^{2 \pi i/5}$.  Together these generate a semi-direct product $S_5 \ltimes (\Z/5)^4$, because the diagonal $\Z/5$ acts trivially.  This is the whole automorphism group by \cite{shioda} or \cite[Lem.~3.2]{balazs}.  We begin with some basic observations about its action on $H^3(Q_0,\Z)$.

\begin{lem}
For $\vec\imath = (i_0, \dotsc, i_4) \in (\Z/5)^5$, let $\phi_{\vec\imath}$ be the automorphism of $Q_0$ defined by
\[ \phi_{\vec\imath}(z_0 : \dotsb : z_4) = (\xi^{i_0} z_0 : \dotsb \colon \xi^{i_4} z_4). \]
For a permutation $\sigma \in S_5$, use the same letter $\sigma$ to denote the corresponding automorphism of $Q_0$.
\begin{enumerate}
\item For $\vec\imath, \vec\jmath \in (\Z/5)^5$ we have $\langle \phi_{\vec\imath}^* \alpha,\, \phi_{\vec\jmath}^* \alpha \rangle = \langle \phi_{\vec\imath - \vec\jmath}^*\, \alpha,\, \alpha \rangle$.

\item For $\sigma \in S_5$ we have $\sigma^* \alpha = \pm \alpha$ depending on whether $\sigma$ is even or odd.

\item For $\vec\imath \in (\Z/5)^5$ and $\sigma \in S_5$ we have $\langle \phi_{\vec\imath}^* \alpha,\, \alpha \rangle = \langle \phi_{\sigma(\vec\imath)}^* \alpha,\, \alpha \rangle$.
\end{enumerate}
\end{lem}
\begin{proof}
(a) This is because $\phi_{\vec\jmath}^{-1} = \phi_{-\vec\jmath}$ preserves the intersection pairing. \medskip

(b) We see that $\sigma$ takes $Q_0(\R)$ into itself, so $\sigma^* \alpha = \pm \alpha$ depending on whether $\sigma$ preserves or reverses the orientation of $Q_0(\R)$.  This determines a group homomorphism from $S_5$ to $\{ \pm 1 \}$, so it is either the sign of the permutation or always $+1$, and to determine which it is we can just look at one transposition.  Let $\sigma$ be the transposition that exchanges the first two coordinates, and consider its action on the tangent space to $Q_0(\R)$ at the fixed point $(0:0:0:-1:1)$, working in the affine patch $z_4 = 1$.  Taking the gradient of the defining equation, we find that this tangent space is $(*,*,*,0) \subset \R^4$, and that $\sigma$ reverses the orientation. \medskip

(c) We find that
\[ \sigma^{-1} \circ \phi_{\vec\imath} \circ \sigma = \phi_{\sigma(\vec\imath)}, \]
or equivalently,
\[ \phi_{\vec\imath} \circ \sigma = \sigma \circ \phi_{\sigma(\vec\imath)}.\]
Thus on cohomology, we have
\[ \sigma^* \circ \phi_{\vec\imath}^* = \phi_{\sigma(\vec\imath)}^* \circ \sigma^*. \]
Now
\[ \langle \phi_{\vec\imath}^* \alpha ,\, \alpha \rangle
= \langle \sigma^* \phi_{\vec\imath}^* \alpha ,\, \sigma^* \alpha \rangle
= \langle \phi_{\sigma(\vec\imath)}^* \sigma^* \alpha ,\, \sigma^* \alpha \rangle
= (\pm 1)^2 \langle \phi_{\sigma(\vec\imath)}^* \alpha ,\, \alpha \rangle. \qedhere \]
\end{proof}

\subsection{Topology of the intersections}

Our next step toward computing $\langle \phi_{\vec\imath}^* \alpha,\, \alpha \rangle = \langle \alpha,\, \phi_{\vec\imath*} \alpha \rangle$ is to analyze the intersection of $Q_0(\R)$ with its translate $\phi_{\vec\imath}(Q_0(\R))$.  This section and the next agree with \cite[\S2.3.2]{bdlr} and \cite[\S3.4]{mw}.

\begin{prop} \label{partition_prop}
The intersection $Q_0(\mathbb R) \cap \phi_{\vec\imath}(Q_0(\R))$ is the same as the intersection of $Q_0(\R)$ with the fixed locus of the action of $\phi_{\vec\imath}$ on $\C\P^4$.  Its topology is determined by the dimensions of the eigenspaces of the corresponding linear action on $\C^5$, or equivalently by the repetitions of entries in $\vec\imath$, as follows:
\pagebreak 
\[ \begin{array}{c|c}
1+1+1+1+1 & \varnothing \\
2+1+1+1 & \text{point} \\
2+2+1 & \text{two points} \\
3+1+1 & S^1 \\
3+2 & S^1 \cup \text{point} \\
4+1 & \R\P^2 \\
5 & \R\P^3
\end{array} \]
\end{prop}
\begin{proof}
Let $\tau\colon Q_0 \to Q_0$ be the anti-holomorphic involution given by complex conjugation:
\[ \tau(z_0 : \dotsb : z_4) = (\bar z_0 : \dotsb : \bar z_4). \]
Thus $Q_0(\R)$ is the fixed locus of $\tau$, and we find that $\phi_{\vec\imath}(Q_0(\R))$ is the fixed locus of $\phi_{\vec\imath} \circ \tau \circ \phi_{\vec\imath}^{-1}$.  Moreover we find that
\[ \tau \circ \phi_{\vec\imath}^{-1} = \tau \circ \phi_{-\vec\imath} = \phi_{\vec\imath} \circ \tau, \]
so $\phi_{\vec\imath}(Q_0(\R))$ is the fixed locus of $\phi_{2\vec\imath} \circ \tau$.

Next we observe that a point in the fixed locus of $\tau$ lies in the fixed locus of $\phi_{2\vec\imath} \circ \tau$ if and only if it lies in the fixed locus of $\phi_{2\vec\imath}$, which is the same as the fixed locus of $\phi_{\vec\imath}$ because the order of $\phi_{\vec\imath}$ is odd.  The fixed locus of a linear automorphism of $\C\P^4$ is the union of the projectivizations of the eigenspaces of the corresponding automorphism of $\C^5$, and because $\phi_{\vec\imath}$ acts diagonally, these are coordinate subspaces, with the dimensions depending on the repetitions in $\vec\imath$.

So for example, if $\vec\imath = (2,2,3,3,3)$, then the fixed locus of
\[ \phi_{\vec\imath}(z_0:z_1:z_2:z_3:z_4) = (\xi^2 z_0 : \xi^2 z_1 : \xi^3 z_2 : \xi^3 z_3 : \xi^3 z_4) \]
has two connected components: the line $(*:*:0:0:0)$, which meets $Q_0(\R)$ in the single point $(1:-1:0:0:0)$, and the plane $(0:0:*:*:*)$, which meets $Q_0(\R)$ in a Fermat quintic curve, homeomorphic to $\R\P^1$ or $S^1$.  Or if $\vec\imath = (2,3,3,3,3)$, then the fixed locus of $\phi_{\vec\imath}$ consists of the point $(1:0:0:0:0)$, which is not in $Q_0$, and the hyperplane $(0:*:*:*:*)$, which meets $Q_0(\R)$ in a Fermat quintic surface, homeomorphic to $\R\P^2$.  The other partition types are similar.
\end{proof}

\subsection{Intersection numbers}

Having determined the intersections $Q_0(\R) \cap \phi_{\vec\imath}(Q_0(\R))$, we study the cup product of the corresponding classes in $H^3(Q_0,\Z)$.  We will apply the following proposition to $M = Q_0$ with the Fubini--Study metric; $L = Q_0(\R)$, which is Lagrangian because complex conjugation flips the sign of the Fubini--Study K\"ahler form; and $\phi = \phi_{\vec\imath}$ for various $\vec\imath$, which preserve the Fubini--Study metric because they come from unitary matrices.  In that setting, the imaginary part of the derivative $d\phi_p$ is just what one would guess.

\begin{prop} \label{orientation_prop}
Let $(M^{2n}, \omega, g, J)$ be a compact K\"ahler manifold, with its standard orientation from the volume form $\omega^n$.  Let $L \subset M$ be a closed, oriented Lagrangian submanifold, and let $\lambda \in H^n(M,\Z)$ be its Poincar\'e dual.
\begin{enumerate}
\item The self-intersection number
\[ \langle \lambda,\, \lambda \rangle  = (-1)^{\binom{n}{2}} \cdot \chi(L), \]
where $\chi$ is the topological Euler characteristic.\footnote{The same fact appears in \cite[Lem.~2.1]{zhang}.}

\item Let $\phi\colon M \to M$ be an automorphism preserving $\omega$, $g$, and $J$.  Suppose that $L$ and $\phi(L)$ intersect transversely, and that all the intersection points are fixed points of $\phi$.  Then the intersection number
\[ \langle \phi^* \lambda,\, \lambda \rangle  = \sum_{p \in L \cap \phi(L)} (-1)^{\binom{n}{2}} \sign \det \Im d\phi_p, \]
where $d\phi_p\colon T_p M \to T_p M$ is the derivative, and its imaginary part is taken with respect to the identification $T_p L \otimes \C \cong T_p M$: precisely, if we use the isomorphism $T_p L \oplus T_p L \to T_p M$ given by $(v,w) \mapsto v + Jw$ to write $d\phi_p$ in the block form $\left( \begin{smallmatrix} A & -B \\ B & A \end{smallmatrix} \right)$, then $\Im d\phi_p$ means $B$.

\item More generally, suppose that $L$ and $\phi(L)$ intersect cleanly, meaning that each connected component $F$ of $L \cap \phi(L)$ is a submanifold, and $T_p F = T_p L \cap T_p \phi(L)$ for all $p \in F$.  Continue to assume that all the intersection points are fixed points of $\phi$.  Then the intersection number
\[ \langle \phi^* \lambda,\, \lambda \rangle = \sum_{\substack{\mathrm{connected} \\ \mathrm{components} \\ F \subset L \cap \phi(L)}} (-1)^{f(n-f) + \binom{n}{2}} \cdot \sign \det \Im d\phi_p|_{(T_p F \otimes \C)^\perp} \cdot \chi(F), \]
where $f = \dim F$, and the imaginary part of $d\phi_p\colon T_p M \to T_p M$ is taken with respect to the identification $T_p L \otimes \C \cong T_p M$ as before, but now we restrict it to the $g$-orthogonal of $T_p F \otimes \C \subset T_p M$, on which it is invertible.
\end{enumerate}
\end{prop}

\begin{proof}
We remark that (a) and (b) are special cases of (c), but it is clearer to prove them one at a time. \bigskip

(a) The intersection number $\langle \lambda,\, \lambda \rangle$ is found by integrating the Euler class of the normal bundle $NL$ over $L$.  Because $L$ is Lagrangian, we find that $NL = J(TL)$: at each point $p \in L$, the $\omega$-orthogonal to $T_p L$ in $T_p M$ is $T_p L$, so the $g$-orthogonal to $T_p L$ is $J(T_p L)$.  It remains to determine whether the orientation of $NL$ determined by the splitting $TM = TL \oplus NL$ agrees with the one given by the isomorphism $J\colon TL \to NL$.

Let $p \in L$, and let $v_1, \dotsc, v_n$ be an oriented basis for $T_p L$.  Our two orientations of $NL$ will agree if and only if the standard orientation of $T_p M$, which is determined by the basis
\begin{equation} \label{std_basis}
v_1, J(v_1),\ \dotsc,\ v_n, J(v_n),
\end{equation}
agrees with the one determined by the basis
\begin{equation} \label{alt_basis}
v_1, \dotsc, v_n,\ J(v_1), \dotsc, J(v_n).
\end{equation}
We see that the number of transpositions required to turn \eqref{std_basis} into \eqref{alt_basis} is $1 + 2 + \dotsb + (n-1) = \binom{n}{2}$. \bigskip

(b) Let $p \in L \cap \phi(L)$.  We have assumed that $\phi(p) = p$.  Let $v_1, \dotsc, v_n$ be an oriented basis for $T_p L$, so $d\phi_p(v_1), \dotsc, d\phi_p(v_n)$ is an oriented basis for $T_p \phi(L)$.  The contribution of $p$ to the intersection number
\[ \langle \phi^* \lambda,\, \lambda \rangle = \langle \lambda,\, \phi_*\lambda \rangle \]
is $+1$ if the orientation of $T_p M$ determined by the basis
\begin{equation} \label{third_basis}
v_1, \dotsc, v_n,\ d\phi_p(v_1), \dotsc, d\phi_p(v_n)
\end{equation}
agrees with the standard orientation determined by \eqref{std_basis} above, and $-1$ if it disagrees.  The change of basis from \eqref{third_basis} to \eqref{alt_basis} is given by the block matrix $\left( \begin{smallmatrix} 1 & A \\ 0 & B \end{smallmatrix} \right)$, where $B$ is as in the statement of the proposition, so it changes the orientation by $\sign \det B$.  The change of basis from \eqref{alt_basis} to \eqref{std_basis} again changes the orientation by $(-1)^{\binom{n}{2}}$. \bigskip

(c) Let $F$ be a connected component of $L \cap \phi(L)$, and consider the ``excess normal bundle''
\[ E = \frac{TM|_F}{TL|_F + T\phi(L)|_F}. \]
If $L$ and $\phi(L)$ intersected transversely then $E$ would be zero, but since they intersect cleanly, its rank is $f = \dim F$.  If $F$ is orientable, then an orientation of $F$ determines an orientation of $E$ by writing
\begin{equation} \label{orientation_of_E}
TM|_F = TF \oplus \frac{TL|_F}{TF} \oplus \frac{T\phi(L)|_F}{TF} \oplus E.
\end{equation}
The contribution of $F$ to the intersection number
\[ \langle \phi^* \lambda,\, \lambda \rangle = \langle \lambda,\, \phi_*\lambda \rangle \]
is given by integrating the Euler class of $E$ over $F$; this fact, which ought to be standard, is given at various levels of generality in \cite[Thm.~VII]{bott}, \linebreak 
\cite[Prop.~3.3]{quillen}, \cite[Prop.~2.2]{ronga}, and \cite[Prop.~3.3]{rs}.  Note that if we reverse the orientation of $F$ then we also reverse the induced orientation of the excess bundle $E$, so the answer is well-defined.

Because $L$ and $\phi(L)$ are Lagrangian and intersect cleanly, we find that $E$ is isomorphic to $J(TF)$: at each point $p \in F$, the $\omega$-orthogonal to $T_p L + T_p \phi(L)$ in $T_p M$ is the intersection of the $\omega$-orthogonals to $T_p L$ and $T_p \phi(L)$, hence is $T_p L \cap T_p \phi(L) = T_p F$, so the $g$-orthogonal to $T_p L + T_p \phi(L)$ is $J(T_p F)$.  It remains to determine whether the orientation of $E$ determined by the splitting \eqref{orientation_of_E} agrees with the one given by the isomorphism $J\colon TF \to E$, and to deal with the case where $F$ is not orientable.

Let $p \in F$.  We have assumed that $\phi$ acts on $F$ as the identity.  Let $v_1, \dotsc, v_f$ be an oriented basis for $T_p F$, and extend it with $v_{f+1}, \dotsc, v_n$ to get an oriented basis for $T_p L$.  Our two orientations of $E$ will agree if and only if the orientation of $T_p M$ determined by the basis
\begin{equation} \label{fourth_basis}
v_1, \dotsc, v_f,\ 
v_{f+1}, \dotsc, v_n,\ 
d\phi_p(v_{f+1}), \dotsc, d\phi_p(v_n),\ 
J(v_1), \dotsc, J(v_f)
\end{equation}
agrees with the standard orientation determined by \eqref{std_basis} above.  We can get from \eqref{fourth_basis} to
\[ v_1, \dotsc, v_f, v_{f+1}, \dotsc, v_n, J(v_1), \dotsc, J(v_f), d\phi_p(v_{f+1}), \dotsc, d\phi_p(v_n) \]
by $f(n-f)$ transpositions.  The change of orientation from there to \eqref{alt_basis} relates to the imaginary part of $d\phi_p$ as before, but now only through its action on $v_{f+1}, \dotsc, v_n$.  And the change of basis from \eqref{alt_basis} to \eqref{std_basis} again contributes $(-1)^{\binom{n}{2}}$.

Finally, if $F$ is not orientable, we can choose a tubular neighborhood of $F$, work on the double cover where the preimage of $F$ becomes orientable, and then divide by 2.
\end{proof}

\begin{rmk}
In part (c), the sign $(-1)^{f(n-f) + \binom{n}{2}}$ could also have been written as $(-1)^{\binom{f}{2} + \binom{n-f}{2}}$.  If, rather than fixing $F$ pointwise, $\phi$ only maps $F$ into itself and fixes some point $p \in F$, then we can still do business, but the intersection number gets a further factor of $\sign \det d\phi_p|_{T_p F}$.
\end{rmk}

Returning now to our calculation of $\langle \phi_{\vec\imath}^* \alpha,\, \alpha \rangle$ and referring to the possible intersections listed in Proposition \ref{partition_prop}, we note that $\chi(S^1) = \chi(\R\P^3) = 0$, so we only need to find the intersection numbers when we are dealing with an isolated point of intersection or with an $\R\P^2$.

For an isolated point of intersection, it is enough to consider
\[ \langle \phi_{(a,b,c,0,0)}^* \alpha,\, \alpha \rangle \]
where $a, b, c \ne 0$, and the fixed point $p = (0:0:0:-1:1)$.  Working in the affine patch $z_4 = 1$ and taking the gradient of the defining equation, we find that the complex tangent space $T_p(Q_0)$ is $(*,*,*,0) \subset \C^4$, and $d\phi_{(a,b,c,0,0)}$ acts with eigenvalues $\xi^a$, $\xi^b$, and $\xi^c$.  Thus by Proposition \ref{orientation_prop}(b), the contribution of $p$ to the intersection number is $(-1)^{3+\nu}$, where $\nu$ counts how many of $a$, $b$, and $c$ are equal to 3 or 4.

For the case where the intersection is $\R\P^2$, it is enough to consider
\[ \langle \phi_{(a,0,0,0,0)}^* \alpha,\, \alpha \rangle \]
where $a \ne 0$, and the same fixed point $p = (0:0:0:-1:1)$.  Now we find that the eigenvalues are $\xi^a$, 1, and 1, so  by Proposition \ref{orientation_prop}(c) the intersection number is $-1$ if $a = 1$ or 2 and $+1$ if $a = 3$ or 4.

\subsection{Matrix computations} \label{matrices}

Take all the classes $\phi_{\vec\imath}^* \alpha$, with the last component of $\vec\imath$ equal to zero on account of the redundancy $\phi_{\vec\imath} = \phi_{\vec\imath + (1,1,1,1,1)}$.  This gives 625 elements of $H^3(Q_0,\Z)$, and thus a homomorphism
\[ A\colon \Z^{625} \to H^3(Q_0,\Z) = \Z^{204}. \]
Theorem \ref{real_parts} asserts that $A$ is surjective, which we will now prove.

Choose a basis of $H^3(Q_0,\Z)$, and let $P$ be the Gram matrix of the intersection pairing.  It is a skew-symmetric matrix with $\lvert \det P \rvert = 1$ by Poincar\'e duality, and in fact $\det P = 1$ because the determinant of a skew-symmetric matrix is the square of its Pfaffian.  If $P'$ is the $625 \times 625$ matrix of intersection pairings between the $\phi_{\vec\imath}^* \alpha$ that we compute by the considerations of the previous section, then
\[ P' = A^\top \circ P \circ A. \]

\begin{lem}
The map $A$ is surjective if and only if $\rank P' = 204$ and $\coker P'$ is torsion-free.
\end{lem}
\begin{proof}
We display the diagram
\[ \xymatrix{
\Z^{625} \ar[r]_-A \ar@/^2em/[rrr]^{P'} & H^3(Q_0,\Z) \ar[r]_-P^-\cong & H^3(Q_0,\Z)^* \ar[r]_-{A^\top} & \Z^{625*}.
} \]
If $A$ is surjective then $A^\top$ is injective, $\coker A^\top$ is torsion-free, $\rank P' = 204$, and $\coker P' = \coker A^\top$.  For the converse, suppose that $\rank P' = 204$.  Then $\rank A = \rank A^\top = 204$, $\coker A$ is finite, and $A^\top$ is injective.  The cokernel of $A^\top$ is isomorphic to $\Z^{421}$ plus a torsion group isomorphic to $\coker A$, and using the fact that $P$ is an isomorphism, we find that $\coker P'$ is an extension of $\coker A^\top$ by $\coker A$; thus $A$ is surjective if and only if $\coker P'$ is torsion-free.
\end{proof}

The ancillary file \verb|quintic.magma| includes Magma code to build the matrix $P'$ and verify that its Smith normal form is a diagonal matrix with 204 ones and 421 zeroes down the diagonal.  Thus $A$ is surjective by the preceding lemma, and Theorem \ref{real_parts} is proved.  We also check that no two rows of $P'$ are equal, so our 625 classes $\phi_{\vec\imath}^* \alpha$ are pairwise distinct.  \medskip

For our application we will need explicit matrices $A$ and $P$ satisfying $P' = A^\top \circ P \circ A$, and \verb|quintic.magma| includes code to find them as follows.  First, find a matrix $B\colon \Z^{204} \to \Z^{625}$ with the same image as $P'$.  Next, find a matrix $A\colon \Z^{625} \to \Z^{204}$ that solves the equation $P' = B \circ A$.  Last, find a matrix $P\colon \Z^{204} \to \Z^{204}$ that solves the equation $B = A^\top \circ P$.  As a sanity check, we verify that $P$ is skew symmetric, and that $\det P = 1$.  The entries of $A$ and $B$ range from $-2$ to 2, while those of $P$ range from $-48$ to 48.

We chose Magma \cite{magma} because it has fast algorithms for dealing with large integer matrices; the disadvantage is that Magma works with row vectors rather than column vectors, viewing the image of a matrix as the row space rather than the column space, and the composition $B \circ A$ as $A * B$ rather than $B * A$, which is confusing for someone who's used to the opposite convention.\medskip

We will also need matrices representing the action of various automorphisms of $Q_0$ on $H^3(Q_0, \Z)$, which we obtain as follows.  For $\phi_{\vec\imath}^*$, we first get a $625 \times 625$ permutation matrix $\Phi'_{\vec\imath}$ that represents $\phi_{\vec\jmath}^* \alpha \mapsto \phi_{\vec\imath}^* \phi_{\vec\jmath}^* \alpha$.  Then we seek a $204 \times 204$ matrix $\Phi_{\vec\imath}$ that satisfies $\Phi_{\vec\imath} \circ A = A \circ \Phi'_{\vec\imath}$.  In fact we have to ask the computer to solve the transposed equation.  The entries in $\Phi_{\vec\imath}$ typically range from $-2$ to 2.  We can deal with $\sigma^*$ for $\sigma \in S_5$ similarly.


\section{The Torelli example}

\subsection{Description of the family} \label{asmosz_setup}

We review the construction of Aspinwall and Morrison \cite{am}, using Szendr\H oi's notation from \cite{balazs}; in particular, we use $t$ as the parameter rather than $\psi$.

Start with the Dwork pencil of quintic 3-folds:
\[ Q_t = \{ z_0^5 + z_1^5 + z_2^5 + z_3^5 + z_4^5 = 5t \cdot z_0 z_1 z_2 z_3 z_4 \} \subset \C\P^4. \]
For $t=0$ this is the Fermat quintic.  For $t \notin \{ 1, \xi, \xi^2, \xi^3, \xi^4 \}$ it is smooth.  As $t \to \infty$ it degenerates to five hyperplanes.

Consider three automorphisms $g_1, g_2, g_3 \in \Aut(Q_t)$
defined by
\begin{align}
g_1(z_0:z_1:z_2:z_3:z_4) &= (z_0 : \xi z_1 : \xi^2 z_2 : \xi^3 z_3 : \xi^4 z_4), \notag \\
g_2(z_0:z_1:z_2:z_3:z_4) &= (z_0 : \xi z_1 : \xi^3 z_2 : \xi z_3 : z_4), \label{def_of_G} \\
g_3(z_0:z_1:z_2:z_3:z_4) &= (z_1 : z_2 : z_3 : z_4 : z_0). \notag
\end{align}
Let $H$ be the group generated by $g_1$ and $g_2$, which is isomorphic to $\Z/5 \times \Z/5$; let $G$ be the group generated by $g_1$, $g_2$, and $g_3$, which is isomorphic to the semi-direct product $(\Z/5 \times \Z/5) \rtimes \Z/5$; and let $K = G/H \cong \Z/5$, which is generated by $g_3$.

The action of $H$ preserves the holomorphic volume form on $Q_t$ and is free except at the 50 points, each with stabilizer $\Z/5$, as we discuss in more detail in \S\ref{1/5(1,1,3)}.  Thus the quotient
\[ \bar Z_t := Q_t/H \]
has ten isolated singularities, and Nakamura's $H$-Hilbert scheme provides a crepant resolution \cite[Thm.~0.1]{nakamura}
\[ \bar Z_t \leftarrow Z_t := H\text{-}\!\Hilb(Q_t) \]
which is again a smooth Calabi--Yau 3-fold.  Szendr\H oi showed in \cite[Prop.~1.7]{balazs} that $Z_t$ is simply connected, and we give an easier proof in Proposition \ref{simply_connected} below.

The residual action of $K$ on $\bar Z_t$ and $Z_t$ is free, so the quotient
\[ \bar Y_t := Q_t/G = \bar Z_t/K \]
has two isolated singularities, and its crepant resolution
\[ \bar Y_t \leftarrow Y_t := Z_t/K = G\text{-}\!\Hilb(Q_t) \]
has $\pi_1(Y_t) = \Z/5$.

In \cite[\S2]{balazs}, Szendr\H oi observed that the pullback $H^3(Y_t,\Q) \to H^3(Z_t,\Q)$ is an isomorphism of polarized Hodge structures (up to multiplying the intersection form by 5), and that the map
\begin{equation} \label{g}
g(z_0:z_1:z_2:z_3:z_4) = (\xi^{-1} z_0 : z_1 : z_2 : z_3 : z_4)
\end{equation}
induces an isomorphism $Z_t \xrightarrow\cong Z_{\xi t}$, so we have $H^3(Y_t, \Q) \cong H^3(Y_{\xi^a t},\Q)$ for $a = 1, 2, 3, 4$.  Our goal in this second half of the paper is to prove the main theorem from the introduction: if $t$ is very general, then there are no such isomorphisms with $\Z$ coefficients, even neglecting the polarization.  In fact we will prove that if $t$ is very general, then scalar multiples of $(g^*)^a$ are the \emph{only} Hodge isomorphisms from $H^3(Y_{\xi^a t},\Q)$ to $H^3(Y_t, \Q)$, but they cannot preserve the integral structures.

\subsection{Crepant resolutions and fundamental groups} \label{1/5(1,1,3)}

The action of $H$ on $Q_t$ is free except at the 50 points given by
\[ (1:-\xi^a:0:0:0) \in Q_t \]
and its permutations, where $a$ runs from 0 to 4.  The stabilizer of each point is isomorphic to $\Z/5$: for example, the stabilizer of $(1:-\xi^a:0:0:0)$ is generated by $g_2 g_1^{-1}$, which we write out for later reference:
\[ g_2 g_1^{-1}(z_0:z_1:z_2:z_3:z_4) = (z_0 : z_1 : \xi z_2 : \xi^3 z_3 : \xi z_4). \]

The resulting singularities of $\bar Z_t = Q_t/H$ are of type $\frac{1}{5}(1,1,3)$ in Reid's notation \cite[\S4.3]{ypg}, meaning that they are locally isomorphic (in the analytic or \'etale topology) to $\C^3$ modulo $\Z/5$ acting on the coordinates by $\xi$, $\xi$, and $\xi^3$.  The quickest way to see this is to project onto a coordinate 3-space: for example, at the point $(1:-\xi^a:0:0:0)$, if we work in the affine patch $z_0 = 1$, then the map $(z_1,z_2,z_3,z_4) \mapsto (z_2,z_3,z_4)$ gives a local isomorphism onto a neighborhood of 0 in $\C^3$, intertwining the action of the stabilizer $\langle g_2 g_1^{-1} \rangle$ with the action of $\Z/5$ that scales $z_2$ by $\xi$, $z_3$ by $\xi^3$, and $z_4$ by $\xi$.

Each exceptional divisor of the crepant resolution $Z_t \to \bar Z_t$ consists of a $\P^2$ and a Hirzebruch surface $\F_3$, intersecting along a $\P^1$ that lies on the $\P^2$ as a straight line and on the $\F_3$ as the negative section \cite[Prop.~1.4]{balazs}.  A systematic approach to the relevant toric calculations was given by Craw and Reid in \cite{cr}.

\begin{prop} \label{simply_connected}
$\bar Z_t$ and $Z_t$ are simply connected.
\end{prop}
\begin{proof}
By the Lefschetz hyperplane theorem, $Q_t$ is simply connected, so by a result of Armstrong \cite{armstrong}, the fundamental group of $\bar Z_t = Q_t/H$ is the quotient of $H$ by the normal subgroup generated by the stabilizers of all fixed points.  We find that these stabilizers generate $H$, so $\pi_1(\bar Z_t) = 0$.

To see that $\pi_1(Z_t) = 0$, we can apply a result of Koll\'ar on fundamental groups of quotient singularities \cite[Thm.~7.8.1]{kollar_shafarevich}, or we can make a direct argument using van Kampen's theorem and the fact that the exceptional divisors $\P^2 \cup \F_3$ are simply connected.
\end{proof}

\subsection{Analysis of \texorpdfstring{$H^3(-,\Z)$}{H3(-,Z)}} \label{H3_section}

Turning to the action of $H$ and $G$ on $H^3(Q_t,\Z) \cong \Z^{204}$, we claim that the invariant parts satisfy
\begin{equation} \label{H^3(Q_t)^H}
H^3(Q_t,\Z)^H = H^3(Q_t,\Z)^G \cong \Z^4,
\end{equation}
and in particular that $K$ acts trivially on $H^3(Q_t,\Z)^H$.  This follows from a calculation in the Jacobi ring of $Q_0$ via Griffiths' residue calculus \cite[\S6.1.3]{voisin_book2}; we include Macaulay2 code to carry out the calculation in the ancillary file \verb|residues.m2|.  It can also be checked using the methods of the first half of this paper, and we include code in \verb|quintic.magma|.

In this section we will see that the torsion-free quotient $H^3(Z_t,\Z)_\text{tf}$ is naturally isomorphic to \eqref{H^3(Q_t)^H}, and $H^3(Y_t,\Z)_\text{tf}$ to an index-25 subgroup of \eqref{H^3(Q_t)^H}.

Quotients by free group actions are easier to study than ones with fixed points, so let $Q_t^\circ \subset Q_t$ be the complement of the 50 points where $H$ does not act freely; let $Z_t^\circ \subset Z_t$ be the complement of the ten exceptional divisors, which is isomorphic to the smooth locus of $\bar Z_t$ and to $Q_t^\circ/H$; and let $Y_t^\circ \subset Y_t$ be the complement of the two exceptional divisors, which is isomorphic to the smooth locus of $\bar Y_t$ and to $Q_t^\circ/G$ and $Z_t^\circ/K$.  We arrange the inclusions, quotient maps, and resolutions in a diagram, naming some of them:
\[ \xymatrix{
Q_t^\circ \ar@{^(->}[rr] \ar@{->>}[d]^{\varpi^\circ} & & Q_t \ar@{->>}[d]^\varpi \\
Z_t^\circ \ar@{->>}[d]^{\pi^\circ} \ar@{^(->}[r] & Z_t \ar@{->>}[r] \ar@{->>}[d]^\pi & \bar Z_t \ar@{->>}[d]^{\bar\pi} \\
Y_t^\circ \ar@{^(->}[r] & Y_t \ar@{->>}[r] & \bar Y_t & 
} \]
All the vertical maps, apart from $\varpi$ in the top right, are quotients by free group actions.  The restriction maps $H^i(Q_t,\Z) \to H^i(Q_t^\circ,\Z)$ are isomorphisms for all $i \le 4$.
 
\begin{prop} \label{open_sets}
The restriction maps
\begin{align*}
H^3(Z_t, \Z) &\to H^3(Z_t^\circ, \Z) \\
H^3(Y_t, \Z) &\to H^3(Y_t^\circ, \Z)
\end{align*}
are isomorphisms.
\end{prop}
\begin{proof}
We prove it for $Z_t$; the proof for $Y_t$ is the same.

Let $E_i = \P^2 \cup \F_3$ be one of the exceptional divisors.  The Mayer--Vietoris sequence gives the homology and cohomology of $E_i$:
\[ \begin{array}{cccccccc}
H_*(E_i) & = & \Z & 0 & \Z^2 & 0 & \Z^2, \\[1ex]
H^*(E_i) & = & \Z & 0 & \Z^2 & 0 & \Z^2.
\end{array} \]

Each singular point $p_i \in \bar Z_t$ has a neighborhood $\bar U_i$ (in the analytic topology) isomorphic to an open ball around the origin in $\C^3$ modulo the action of $\Z/5$ that multiplies the coordinates by $\xi$, $\xi$, and $\xi^3$.  Let $U_i$ be the preimage of $\bar U_i$ in $Z_t$, which deformation retracts onto the exceptional divisor $E_i$ lying over $p_i$.  The intersection
\[ Z_t^\circ \cap U_i = U_i \setminus E_i = \bar U_i \setminus p_i \]
deformation retracts onto a lens space $S^5/(\Z/5)$, so a standard calculation with the Cartan--Leray spectral sequence \cite[Thm.~6.10.10]{weibel} gives its homology and cohomology:
\[ \begin{array}{cccccccc}
H_*(Z_t^\circ \cap U_i) & = & \Z & \Z/5 & 0 & \Z/5 & 0 & \Z, \\[1ex]
H^*(Z_t^\circ \cap U_i) & = & \Z & 0 & \Z/5 & 0 & \Z/5 & \Z.
\end{array} \]
We remind the reader that
\[ H^i(\Z/5,\,\Z) = \begin{cases}
\Z & i=0 \\
0 & i = 1,3,5,\dotsc \\
\Z/5 & i = 2,4,6,\dotsc.
\end{cases} \]

The long exact sequence of the pair $(Z_t, Z_t^\circ)$ includes
\[ H^3(Z_t, Z_t^\circ) \to H^3(Z_t) \to H^3(Z_t^\circ) \to H^4(Z_t, Z_t^\circ). \]
By Lefschetz duality, the first term is isomorphic to $\bigoplus_{i=1}^{10} H_3(E_i)$, which is zero, so the middle map is injective.  The last term is isomorphic to $\bigoplus H_2(E_i),$ which is not zero, but we will argue that the last map is zero nonetheless.  By excision, the last term is isomorphic to $\bigoplus H^4(U_i, Z_t^\circ \cap U_i)$, and the diagram
\[ \xymatrix{
H^3(Z_t^\circ) \ar[r] \ar[d] & H^4(Z_t, Z_t^\circ) \ar@{=}[d]\\
\bigoplus H^3(Z_t^\circ \cap U_i) \ar[r] & \bigoplus H^4(U_i, Z_t^\circ \cap U_i)
} \]
commutes.  But we have seen that $H^3(Z_t^\circ \cap U_i) = 0$.
\end{proof}

\begin{prop} \label{H^3(Z_t)}
The pullback map
\begin{equation} \label{varpi^circ^*}
\varpi^{\circ*}\colon H^3(Z_t^\circ, \Z) \to H^3(Q_t^\circ, \Z)^H = \Z^4
\end{equation}
is surjective, and its kernel is $\Z/5$.  Thus
\[ H^3(Z_t,\Z) \cong \Z^4 \oplus \Z/5. \]
\end{prop}
\begin{proof}
First we prove that \eqref{varpi^circ^*} is surjective.  Consider the pushforward or norm map
\[ \varpi^\circ_*\colon H^3(Q^\circ_t,\Z) \to H^3(Z^\circ_t,\Z), \]
and observe that the composition
\[ \varpi^{\circ*} \circ \varpi^\circ_*: H^3(Q^\circ_t,\Z) \to H^3(Q^\circ_t,\Z)^H \]
is given by the sum
\[ \sum_{h \in H} h^* = \sum_{i,j=0}^4 (g_1^*)^i \circ (g_2^*)^j. \]
Using the spanning set for $H^3(Q_0,\Z)$ from the first half of the paper, we find that this is surjective for $t=0$ with code included in \verb|quintic.magma|.  Hence it is surjective for all $t$, so $\varpi^{\circ*}$ is surjective.

For the kernel of \eqref{varpi^circ^*}, we use the Cartan--Leray spectral sequence
\[ E_2^{p,q} = H^p(H, H^q(Q_t^\circ)) \Rightarrow H^{p+q}(Z_t^\circ). \]
The $E_2$ page reads as follows:
\begin{equation} \label{E_2_page}
\begin{split}
\xymatrix@=1ex{
\vdots & \vdots & \vdots & \vdots & \vdots & \vdots  \\
3 & H^3(Q_t^\circ)^H & ? & ? & ? & ? & \cdots \\
2 & \Z & 0 & (\Z/5)^2 & \Z/5 & (\Z/5)^3 & \cdots \\
1 & 0 & 0 & 0 & 0 & 0 & \cdots \\
q=0 & \Z & 0 & (\Z/5)^2 & \Z/5 & (\Z/5)^3 & \cdots \\
\ar@{-}[]+<1.6em,.8em>;[rrrrrr]+<1.6em,.8em> 
\ar@{-}[]+<1.6em,.8em>;[uuuuu]+<1.6em,.8em> 
& p=0 & 1 & 2 & 3 & 4 & \cdots
} \end{split}
\end{equation}
In the rows $q=0$ and $q=2$, we have used the fact that $H$ acts trivially on
\[ H^0(Q_t^\circ) = H^0(Q_t) = \Z \qquad \text{and} \qquad H^2(Q_t^\circ) = H^2(Q_t) = \Z, \]
and have computed the cohomology of $H = \Z/5 \times \Z/5$, or of its classifying space $BH = B(\Z/5) \times B(\Z/5)$, from that of $\Z/5$ or $B(\Z/5)$ using the K\"unneth formula.

Our claim about the kernel of \eqref{varpi^circ^*} amounts to saying that the $\Z/5$ at $p=3$, $q=0$ survives to the $E_\infty$ page, or that the differential $d_3: E_3^{0,2} \to E_3^{3,0}$ vanishes.  This differential goes from $H^2(Q_t^\circ)^H$ to $H^3(H, \Z)$, and measures the difference between $H$-invariant line bundles and those that admit an $H$-equivariant structure, hence descend to the quotient.  Now $H^2(Q_t^\circ)$ is generated by the restriction of $\O_{\P^4}(1)$, which admits an $H$-equivariant structure because the action of $H$ on $\P^4$ comes from a linear representation on $\C^5$.
\end{proof}

\begin{prop} \label{H^3(Y_t)}
The image of the pullback map
\begin{equation} \label{pi^circ^*}
\varpi^{\circ*} \circ \pi^{\circ*}\colon H^3(Y_t^\circ, \Z) \to H^3(Q_t^\circ, \Z)^G = H^3(Q_t^\circ, \Z)^H = \Z^4 
\end{equation}
has index 25, and the kernel is $\Z/5$.  Thus
\[ H^3(Y_t,\Z) \cong \Z^4 \oplus \Z/5. \]
\end{prop}
\begin{proof}
We use the Cartan--Leray spectral sequence
\[ E_2^{p,q} = H^p(G, H^q(Q_t^\circ, \Z)) \Rightarrow H^{p+q}(Y_t^\circ,\Z). \]
The $E_2$ page is similar to \eqref{E_2_page}, but the rows $q=2$ and $q=0$ are
\[ \begin{array}{cccccccc}
H^*(G,\Z) & = & \Z & 0 & (\Z/5)^2 & (\Z/5)^2 & (\Z/5)^4 & \cdots
\end{array} \]
which we can get from Lewis \cite[Thm.~6.26]{lewis} or with GAP \cite{GAP4}.  (It is helpful to know that $G$ is the unique non-Abelian group of order 125 and exponent 5.)  Now we claim that the differential on the $E_3$ page from $H^2(Q_t^\circ)^G = \Z$ to $H^3(G,\Z) = (\Z/5)^2$ is not zero.  The restriction of $\O_{\P^4}(1)$ is $G$-invariant but not $G$-equivariant, because the projective representation of $G$ given in \eqref{def_of_G} does not lift to a linear representation, but gives a non-trivial element of the Schur multiplier: for example, we see that while $g_1 \circ g_3$ and $g_3 \circ g_1$ agree as automorphisms of $\P^4$, they differ by a multiple of $\xi$ as automorphisms of $\C^5$.  Thus the image of this differential is $\Z/5$, so its cokernel is $\Z/5$, which survives to the $E_\infty$ page.  Thus the kernel of \eqref{pi^circ^*} is $\Z/5$, and from the diagonal $p+q=4$ we see that its image has finite index.

To show that the index is 25, by Propositions \ref{open_sets} and \ref{H^3(Z_t)} it is enough to show that the image of the pullback map
\[ \pi^*\colon H^3(Y_t,\Z)_\text{tf} \to H^3(Z_t,\Z)_\text{tf} \]
has index 25.  Consider the pushforward or norm map
\[ \pi_*\colon H^3(Z_t,\Z)_\text{tf} \to H^3(Y_t,\Z)_\text{tf}. \]
Because $\pi$ is a 5-sheeted covering space, the composition $\pi_* \circ \pi^*$ is multiplication by $5$.  If we choose bases for both groups and let $A$ be the $4 \times 4$ matrix representing $\pi^*$, then the index of the image is $\lvert \det A \rvert$.  For $\alpha, \beta \in H^3(Y_t,\Z)_\text{tf}$ we have
\[ \langle \pi^* \alpha,\, \pi^* \beta \rangle_{Z_t} = \langle \alpha,\, \pi_* \pi^* \beta \rangle_{Y_t} = 5 \cdot \langle \alpha,\, \beta \rangle_{Y_t}, \]
so if $P_Z$ is the Gram matrix of the intersection pairing on $H^3(Z_t,\Z)_\text{tf}$ in the chosen basis, and $P_Y$ is the Gram matrix for $H^3(Y_t,\Z)_\text{tf}$, then
\[ A^\top P_Z A = 5\cdot P_Y. \]
By Poincar\'e duality we have $\lvert \det P_Z \rvert = \lvert \det P_Y \rvert = 1$, so taking determinants we find that $\lvert \det A \rvert^2 = 5^4$.
\end{proof}

\begin{rmk}
We have avoided using the Cartan--Leray spectral sequence
\[ E_2^{p,q} = H^p(K, H^q(Z_t, \Z)) \Rightarrow H^{p+q}(Y_t,\Z), \]
which Szendr\H oi uses in the proof of \cite[Lem.~2.1]{balazs}, for several reasons.  First, a careful analysis of $H^2(Z_t,\Z) = \Z^{21}$ is needed to decide whether its $K$-module structure is $\Z \oplus (\Z K)^4$, or $\Z^5 \oplus (\Z K / \Z)^4$, or something in between.  Second, while $K$ acts trivially on the torsion subgroup of $H^3(Z_t, \Z) = \Z^4 \oplus \Z/5$ and on its torsion-free quotient, this doesn't imply that it acts trivially on the whole group.  Third, the differentials coming out of $p=0$, $q=3$ on the $E_2$ page and especially the $E_4$ page are hard to determine, as Szendr\H oi remarks in \cite[Rmk.~2.4]{balazs}.  We can deduce some of the answers to these questions from Proposition \ref{H^3(Y_t)}, but we do not know how to get them directly.
\end{rmk}

\subsection{Analysis of Hodge endomorphisms} \label{Hodge_endos}

Take the map $g\colon Q_t \to Q_{\xi t}$ defined in \eqref{g}, and set $t = 0$ to get an automorphism of $Q_0$.  It commutes with $H = \langle g_1, g_2 \rangle$, hence acts on $\bar Z_0$ and $Z_0$, but it does not normalize $G = \langle g_1, g_2, g_3 \rangle$, so it does not descend to $\bar Y_0$ and $Y_0$.

\begin{prop} \label{eigenvalues}
The action of $g^*$ on $H^3(Z_0, \C)$ has eigenvalues $\xi^{-1}$, $\xi^{-2}$, $\xi^{-3}$, and $\xi^{-4}$, with eigenspaces $H^{3,0}(Z_0)$, $H^{2,1}(Z_0)$, $H^{1,2}(Z_0)$, and $H^{0,3}(Z_0)$, respectively.
\end{prop}
\begin{proof}
By Propositions \ref{open_sets} and \ref{H^3(Z_t)} we see that $H^{p,q}(Z_0)$ is identified with the $H$-invariant subspace of $H^{p,q}(Q_0)$ for $p+q=3$.  We have analyzed the latter using Griffiths' residue calculus in \verb|residues.m2|, finding that the $H$-invariant subspaces of $H^{3,0}(Q_0)$, $H^{2,1}(Q_0)$, $H^{1,2}(Q_0)$, and $H^{0,3}(Q_0)$ are generated by the residues of
\begin{align*}
\frac{\Omega}{f} &&
z_0 z_1 z_2 z_3 z_4 \cdot \frac{\Omega}{f^2} &&
(z_0 z_1 z_2 z_3 z_4)^2 \cdot \frac{\Omega}{f^3} &&
(z_0 z_1 z_2 z_3 z_4)^3 \cdot \frac{\Omega}{f^4},
\end{align*}
respectively, where $f = z_0^5 + \dotsb + z_4^5$ and
\[ \Omega = z_0\, dz_1 \wedge dz_2 \wedge dz_3 \wedge dz_4 - z_1 dz_0 \wedge dz_2 \wedge dz_3 \wedge dz_4 + \dotsb. \]
Because $g$ multiplies $z_0$ by $\xi^{-1}$ and leaves the other coordinates alone, we see that it multiplies the meromorphic 4-forms above by $\xi^{-1}$, $\xi^{-2}$, $\xi^{-3}$, and $\xi^{-4}$ respectively.
\end{proof}

\begin{prop} \label{endomorphisms} \ 
\begin{enumerate}
\item The Hodge endomorphism algebra $\End_{\Hdg}(H^3(Z_0, \Q))$ is isomorphic to the cyclotomic field $\Q[\xi]$.
\item For all but countably many values of $t$, $\End_{\Hdg}(H^3(Z_t,\Q)) = \Q$.
\end{enumerate}
\end{prop}
\begin{proof}
(a) On the one hand, $g^*$ is in $\End_{\Hdg}(H^3(Z_0, \Q))$, and from the eigenvalue calculation in Proposition \ref{eigenvalues} we see that the minimal polynomial of $g^*$ is $x^4 + x^3 + x^2 + x + 1$, so it generates a subalgebra isomorphic to $\Q[\xi]$, which we note is 4-dimensional.

On the other hand, $\End_{\Hdg}(H^3(Z_0, \Q))$ is contained in the $(0,0)$-part of $\End_\Q(H^3(Z_0, \Q))$, and because the Hodge numbers of $H^3(Z_0, \Q)$ are $1,1,1,1$, the Hodge numbers of $\End_\Q(H^3(Z_0, \Q))$ are $1,2,3,4,3,2,1$, so $\End_{\Hdg}(H^3(Z_0, \Q))$ cannot be any bigger.  \bigskip

(b) Let $B = \C \setminus \{1,\xi,\xi^2,\xi^3,\xi^4\}$, which is the base of all our families of 3-folds, and let $\tilde B$ be its universal cover.  Given a class $f \in \End_\Q(H^3(Z_0,\Q))$, let $B_f \subset B$ be the set of all $t$ for which parallel transport along some path from 0 to $t$ takes $f$ to a class of type $(0,0)$.  By \cite[Lem.~5.13]{voisin_book2}, this is the image of an analytic subvariety of $\tilde B$, so it is either countable or all of $B$.  We will argue that if $f$ is not a multiple of the identity then $B_f$ is not all of $B$; then the union of all such $B_f$ is a countable union of countable sets, which we discard.

Clearly if $f$ is not of type $(0,0)$ at $t=0$ then $B_f \varsubsetneq B$.  If $f$ is of type $(0,0)$ at $t=0$ and remains so for all $t$, then the composition
\[ T_0(B) \to H^1(TZ_0) \to \Hom_\C(H^{3,0}(Z_0), H^{2,1}(Z_0)) \]
must take values in the $f$-invariant part of the target.  The first map (the Kodaira--Spencer map) is an isomorphism by a straightforward computation using \cite[Lem.~6.15]{voisin_book2}, and the second is an isomorphism by the local Torelli theorem.  So it is enough to show that if $f$ is not a multiple of the identity, then $f$ acts non-trivially on the target, or equivalently, that it acts on the 1-dimensional spaces $H^{3,0}(Z_0)$ and $H^{2,1}(Z_0)$ with different eigenvalues.

Since $f$ is of type $(0,0)$ at $t=0$, by part (a) we can write
\[ f = x + y g^* + z g^{*2} + w g^{*3} \]
for some $x,y,z,w \in \Q$.  By Proposition \ref{eigenvalues}, $f$ acts on $H^{3,0}$ and $H^{2,1}$ with eigenvalues
\begin{align*}
\lambda_{3,0} &:= x + y \xi^{-1} + z \xi^{-2} + w \xi^{-3} \\
\lambda_{2,1} &:= x + y \xi^{-2} + z \xi^{-4} + w \xi^{-6}.
\end{align*}
If these are equal, then $\lambda_{3,0}$ is fixed by the automorphism of the field $\Q[\xi]$ that takes $\xi$ to $\xi^2$; but this automorphism generates the Galois group of $\Q[\xi]$ over $\Q$, so we must have $\lambda_{3,0} \in \Q$, so $f$ is multiplication by a scalar.
\end{proof}

\subsection{End of the proof} \label{finis_laus_deo}

In the statement of our main theorem, let us take ``very general'' to mean that $t$ satisfies the conclusion of Proposition \ref{endomorphisms}(b).  Suppose that for some such $t$ and some $a \in \{1,2,3,4\}$ we had a (possibly unpolarized) Hodge isomorphism $\eta\colon H^3(Y_{\xi^a t}, \Z)_\text{tf} \to H^3(Y_t, \Z)_\text{tf}$.  Tensoring with $\Q$ would give an isomorphism $H^3(Y_{\xi^a t}, \Q) \to H^3(Y_t, \Q)$, and thus $H^3(Z_{\xi^a t}, \Q) \to H^3(Z_t, \Q)$, which we still call $\eta$.  Then $(g^*)^{-a} \circ \eta$ must be multiplication by some rational number $r \in \Q$, so $\eta = r(g^*)^a$.  Thus the map
\[ r(g^*)^a \colon H^3(Z_{\xi^a t}, \Q) \to H^3(Z_t, \Q) \]
would take $\pi^* H^3(Y_{\xi^a t}, \Z)_\text{tf}$ into $\pi^* H^3(Y_t, \Z)_\text{tf}$.  Because $g^*$ is defined across the whole family, this holds for some $t$ if and only if it holds for all $t$, and in particular for $t=0$: the endomorphism $r(g^*)^a$ of $H^3(Z_0, \Q)$ would take $\pi^* H^3(Y_0, \Z)_\text{tf}$ into itself.  Taking fifth powers and recalling that $(g^*)^5 = 1$, we see that multiplication by $r^5$ would take $\pi^* H^3(Y_0, \Z)_\text{tf}$ into itself, so $r^5 = \pm 1$, so $r = \pm 1$.  Thus the automorphism $(g^*)^a$ of $H^3(Z_0, \Z)_\text{tf}$ would take the index-25 subgroup $\pi^* H^3(Y_0, \Z)_\text{tf}$ into itself; but we will argue that it does not.\bigskip

Recalling Propositions \ref{open_sets}, \ref{H^3(Z_t)}, and \ref{H^3(Y_t)}, let us restrict to $Y_0^\circ$ and $Z_0^\circ$, and pull back to $Q_0^\circ$.  So we must show that the automorphism $(g^*)^a$ of $H^3(Q_0^\circ,\Z)^G$ does not take the index-25 subgroup $\varpi^{\circ*} \pi^{\circ*} H^3(Y_0^\circ, \Z)$ into itself.  Because $a$ is relatively prime to the order of $g^*$, it is enough to prove it for $a=1$.

In the first half of the paper we studied the cohomology class $\alpha = [Q_0(\R)] \in H^3(Q_0,\Z)$.  Now we restrict it to $Q_0^\circ$ and choose a singular cocycle $\tilde\alpha \in Z^3_\text{sing}(Q_0^\circ,\Z) \subset C^3_\text{sing}(Q_0^\circ,\Z)$ that represents it.  There is no difficulty in choosing compatible representatives for the classes $\phi_{\vec\imath}^* \alpha$: we just take $\phi_{\vec\imath}^* \tilde\alpha$, because $\phi_{\vec\imath}$ acts on the spaces $Q_0$ and $Q_0^\circ$ and thus on their singular cochain groups, not just on cohomology.  So we have taken the map $A\colon \Z^{625} \to H^3(Q_0, \Z) = H^3(Q_0^\circ, \Z)$ from \S\ref{matrices} and lifted it to a map $\tilde A\colon \Z^{625} \to Z^3(Q_0^\circ,\Z)$, which is still $G$-equivariant, and is now injective because the submanifolds $\phi_{\vec\imath}(Q_0^\circ(\R))$ are distinct.

Because $Q_0^\circ \to Y_0^\circ$ is a quotient by a free action of $G$, the pullback $\varpi^{\circ*} \circ \pi^{\circ*}$ identifies the cochain group $C^3(Y_0^\circ,\Z)$ with the $G$-invariant cochains $C^3(Q_0^\circ,\Z)^G$ by \cite[Lem.~6.10.2]{weibel},\footnote{This reference actually gives the dual statement, that chains downstairs are identified with coinvariant chains upstairs.} and thus identifies the cocycles $Z^3(Y_0^\circ,\Z)$ with the $G$-invariant cocycles $Z^3(Q_0^\circ,\Z)^G$ (although we could not say the same about coboundaries).  In particular, we can get classes in the image of
\[ \varpi^{\circ*} \circ \pi^{\circ*} \colon H^3(Y_0^\circ, \Z) \to H^3(Q_0^\circ, \Z)^G \]
by intersecting our $\Z^{625} \subset Z^3(Q_0^\circ,\Z)$ with the $G$-invariant subgroup and then mapping to $H^3(Q_0^\circ,\Z)^G$.  This is done in \verb|quintic.magma|: the intersection is a $\Z^9 \subset \Z^{625}$, and its image in $H^3(Q_0^\circ, \Z)^G$ has index 125.  So we only get our hands on a subgroup of $\varpi^{\circ*} \pi^{\circ*} H^3(Y_0^\circ, \Z)$, not the whole thing, but this turns out to be enough:

If $g^*$ preserved $\varpi^{\circ*} \pi^{\circ*} H^3(Y_0^\circ, \Z)$, then this index-125 subgroup and its translates by $g^*$, $g^{*2}$, and $g^{*3}$ would at most span a subgroup of index 25 in $H^3(Q_0^\circ,\Z)^G$.  But with \verb|quintic.magma| we find on the contrary that they span all of $H^3(Q_0^\circ, \Z)^G$.

\bibliographystyle{plain}
\bibliography{asmosz}

\scriptsize
\noindent Nicolas Addington \\
Department of Mathematics \\
University of Oregon \\
Eugene, OR 97403-1222 \\
adding@uoregon.edu \\

\noindent Benjamin Tighe \\
Department of Mathematics \\
Pacific Lutheran University \\
12180 Park Ave.\ S \\
Tacoma, WA 98447 \\
btighe@plu.edu

\end{document}